\documentclass[twoside,12pt]{article}
\newcommand{\Q}{\bigskip\par\noindent}

\usepackage{amssymb}
\newcommand{\goth}{\mathfrak}

\def\oqmn{{\cal O}_q(M_n)}
\def\oqmnc{{\cal O}_q(M_n({\Bbb C}))}
\def\rtx{R_t(\mat{X})}

\newcommand{\ideal}[1]{\langle #1 \rangle}
\def\widebar{\overline}
\def\qdot{q^\bullet}

\renewcommand{\k}{{\Bbb K}}
\newcommand{\GK}{{\rm GKdim}\,}

\renewcommand{\det}{{\rm det}}

\newcommand{\mat}[1]{{\bf #1}}

\newcommand{\proof}{\noindent {\it Proof:$\;$}}
\newcommand{\qed}{\hfill \rule{1.5mm}{1.5mm}}
\newcommand{\QED}{\hfill \rule{1.5mm}{1.5mm}\\}

\newtheorem{subtheorem}{Theorem}[section]
\newtheorem{subproposition}[subtheorem]{Proposition}
\newtheorem{sublemma}[subtheorem]{Lemma}
\newtheorem{subcorollary}[subtheorem]{Corollary}
\newtheorem{subremark}[subtheorem]{Remark}

\newtheorem{theorem}{Theorem}[section]
\newtheorem{proposition}[theorem]{Proposition}
\newtheorem{lemma}[theorem]{Lemma}

\newtheorem{relation}[theorem]{Relation}
\newtheorem{relations}[theorem]{Relations}

\newcommand{\titre}{The maximal order property for quantum determinantal
rings} 

\begin{document}
 
\title{{\vspace{-1.5cm} \titre}}

\author{T H Lenagan \ and \ L Rigal\thanks{This research was partially
supported by grants from the Edinburgh and London Mathematical Societies
and  the
European Science Foundation programme `Noncommutative Geometry'.}
}
\date{}
\maketitle

\begin{abstract}  
We develop a method of reducing the size of quantum minors in the
algebra of quantum matrices $\oqmn$. We use the method to show that
the quantum determinantal factor rings of $\oqmnc$ are maximal orders,
for $q$ an element of ${\mathbb C}$ transcendental over ${\mathbb Q}$. 
\end{abstract} 


\section*{Introduction}
Throughout, $\k$ will denote a base field, $q$ a non-zero element of $\k$,
and $m,n$ positive integers. 
We denote by ${\cal O}_q(M_{m,n})$ the quantization of the ring of
regular functions on $m \times n$ matrices with entries in $\k$; it is
the $\k$-algebra generated by $mn$ indeterminates $X_{ij}$, $1 \le i \le
m$ and $1 \le j \le n$, subject to the following relations:

\begin{equation} \label{definition-of-quantum-matrices}
\begin{array}{rcl}
X_{ij}X_{il} &=& qX_{il}X_{ij}, \\
X_{ij}X_{kj} &=& qX_{kj}X_{ij}, \\
X_{il}X_{kj} &=& X_{kj}X_{il},  \\
X_{ij}X_{kl} - X_{kl}X_{ij} &=& (q-q^{-1})X_{il}X_{kj},
\end{array}
\end{equation}
for $1 \le i < k \le m$ and $1 \le j < l \le n$.\\

It will be convenient to use the following notation: setting
$\mat{X}=(X_{ij})_{1\le i \le m,1 \le j\le n}$ we will denote ${\cal
O}_q(M_{m,n})$ by $\k_{q}[\mat{X}]$.  \\

   Let $A$ be any $\k$-algebra, we say that an $m\times n$ 
matrix $(a_{ij})$ with 
entries in $A$ is a {\bf $q$-quantum matrix} if the map
$X_{ij} \mapsto a_{ij}$ induces a homomorphism of algebras
$\k_{q}[\mat{X}] \longrightarrow A$.  If, in addition, this homomorphism
is a monomorphism, then we say that the matrix $(a_{ij})$
is a {\bf generic $q$-quantum matrix}; clearly, $\mat{X}$ is a
generic $q$-quantum matrix, by definition.  \\

   In the case where $m=n$, the {\bf quantum determinant} of $\mat{X}$ is
defined by 
\[
\det_q =\det_{q} \mat{X}:= \sum_{\sigma \in S_{n}}(-q)^{l(\sigma)} 
X_{1,\sigma(1)} \dots X_{n,\sigma(n)};
\] 
from \cite[Theorem 4.6.1]{PW}, 
we know that $\det_q$ is in the centre of ${\cal
O}_q(M_n):={\cal O}_q(M_{n,n})$.  Clearly, if $I \subseteq {\Bbb
N}_m:=\{1,\dots,m\}$ and $J \subseteq {\Bbb N}_n:=\{1,\dots,n\}$ with
$|I|=|J|=t$, that is, $I$ and $J$ both have $t$ elements,  
then the submatrix obtained from $\mat{X}$ by 
keeping the rows
indexed by elements of $I$ and the columns indexed by elements of $J$ is
a generic $q$-quantum matrix; and so, we can speak of its quantum
determinant.  Such an element is called a 
{\bf $t \times t$ quantum minor} of 
$\mat{X}$ and is denoted by $[I|J]$.  In order to simplify the 
notation, if
$I$ and $J$ are given by the explicit list of their elements:
$I=\{i_{1}<\dots<i_{t}\}$ and $J=\{j_{1}<\dots<j_{t}\}$ we will 
use the notation $[i_{t},\dots,i_{1}|j_{1},\dots,j_{t}]$.  It follows at
once from the above that, 
\begin{equation}
\label{semi-centrality-of-minors} \mbox{ for } I \subseteq {\Bbb N}_m,
\; J \subseteq {\Bbb N}_n, \mbox{ with } |I|=|J|=t, \; i\in I\; \mbox{
and } j \in J,\; \; [I|J]X_{ij}=X_{ij}[I|J].  
\end{equation}

We need to use several identities involving quantum minors.  Many of
these are obtained from \cite{PW}.  However, 
it is worth noting that the conventions used in \cite{PW} are slightly
different from ours.  In order to make the notation fit we must replace
$q$ by $q^{-1}$ each time we use a relation from \cite{PW}. 
For convenience, 
some of the identities we use are collected in an Appendix.  \\

   For any positive integer $t$ such that $t \le \min \{m,n\}$, we
denote by ${\cal I}_{t}(\mat{X})$ the two-sided ideal of
$\k_{q}[\mat{X}]$ generated by the $t \times t$ quantum minors of
$\mat{X}$; such ideals will be referred to as {\bf quantum determinantal
ideals}.  For each such $t$, we define the {\bf quantum determinantal
ring} $R_{t}(\mat{X}):=\k_{q}[\mat{X}]/{\cal I}_{t}(\mat{X})$.  It has
been proved in \cite{GL} that $R_{t}(\mat{X})$ is a domain.  Here, we
are interested in the question as to whether or not $\rtx$ is a maximal
order in its division ring of fractions.  This has already been
established for $R_{2}(\mat{X})$ in \cite{R}.  We will show that the
localised ring $\rtx[x_{1n}^{-1}]$, 
where $x_{1n}:= X_{1n} + {\cal I}_{t}(\mat{X})$, 
is a maximal order for each $t$, and
are then able to deduce that $\rtx$ is a maximal order, in the case that
$\k = {\Bbb C}$ and $q$ is an element of ${\mathbb C}$ transcendental over
${\mathbb Q}$.  The question for
general $\k$ and any $q$ remains open: we conjecture that all quantum
determinantal rings are maximal orders.  The problem is a technical one:
it is necessary to show that a certain ideal is a prime ideal.  The
restriction to $q$  an element of ${\mathbb C}$ transcendental over
${\mathbb Q}$ is because the
relevant ideal is shown to be prime in \cite{Lau} for this case.  We can
answer the question for general $0 \neq q \in \k$ for the case of
$R_{n}(\mat{X}) \; ( = \oqmn/\langle \det_q \rangle) $, when $\mat{X}$
is $n\times n$,  and also for the
case $R_{3}(\mat{X})$ for general $m\times n$.  
A proof of the former case is included in this
paper since it is relatively short, and of independent interest.  The
proof of the latter is not given here, since it is somewhat ad-hoc and
rather long. 



\section{A new set of generators for $\k_{q}[\mat{X}]_{X_{1n}}$}
\label{newsetofgen}

It follows at once from (\ref{definition-of-quantum-matrices}) that
$X_{1n}$ is a normal element of $\k_{q}[\mat{X}]$.  So, $\{X_{1n}^{i},
\; i \in {\Bbb N}\}$ is a right denominator set of $\k_{q}[\mat{X}]$;
the corresponding localisation of $\k_{q}[\mat{X}]$ will be denoted
$\k_{q}[\mat{X}]_{X_{1n}}$.

For the remainder of this section we assume that $\min\{m,n\} \ge 2$. In
$\k_{q}[\mat{X}]_{X_{1n}}$, we consider the following elements, for
$2 \le i \le m$ and $1 \le j \le n-1$:
\begin{equation} \label{definition-of-X-prime}
X_{ij}^{'} = X_{ij}-q^{-1}X_{1j}X_{in}X_{1n}^{-1} = -q^{-1}[i,1|j,n]X_{1n}^{-1}.
\end{equation}

Set 
\[
\widetilde{{\bf X}} = 
\left(
\begin{array}{cccc}
X_{11} & \dots & X_{1,n-1} & X_{1n}\\
X_{21}^{\prime} & \dots & X_{2,n-1}^{\prime} & X_{2n}\\
\vdots & & \vdots & \vdots \\
X_{m1}^{\prime} & \dots & X_{m,n-1}^{\prime} & X_{mn}
\end{array}
\right) \quad \mbox{ and } \quad 
{\bf X}^{\prime} = 
\left(
\begin{array}{ccc}
X_{21}^{\prime} & \dots & X_{2,n-1}^{\prime} \\
\vdots & & \vdots \\
X_{m1}^{\prime} & \dots & X_{m,n-1}^{\prime}  
\end{array}
\right)
\]
with the convention that the index of a row in ${\bf X}^{\prime}$ is
actually the index of the corresponding row in $\widetilde{\bf X}$.  So,
for instance $X_{21}^{\prime} \dots X_{2,n-1}^{\prime}$ is referred to
as the row of index $2$ in ${\bf X}^{\prime}$.  \\

   It is clear that the entries of $\widetilde{\bf X}$ together with
$X_{1n}^{-1}$ form a new set of algebra generators for
$\k_q[\mat{X}]_{X_{1n}}$.  Our next aim is to give a description of
$\k_q[\mat{X}]_{X_{1n}}$ based on this set of generators.  \\

In order to do this, we need to know 
the relations between the entries of the matrix
$\widetilde{\mat{X}}$. This problem is dealt with in 
Lemma \ref{X-prime-is-q-quantum} and Proposition
\ref{relations-in-X-tilda} where we obtain a set of relations. 
We will show later that the relations obtained in these results give 
the complete set of relations between entries of $\widetilde{\mat{X}}$. 

\begin{sublemma} \label{X-prime-is-q-quantum}
The matrix $\mat{X}'$ is a $q$-quantum matrix, 
and all its entries commute with $X_{1n}$.
\end{sublemma}

\proof Let $2 \le i,k \le m$ and $1 \le j,l \le n-1$.  From
(\ref{semi-centrality-of-minors}) we see that $X_{1n}$ commutes with the
$2 \times 2$ quantum minor $[i,1|j,n]$.  So
$X_{ij}^{\prime}X_{1n}=X_{1n}X_{ij}^{\prime}$ follows from the second
equality of (\ref{definition-of-X-prime}).  Moreover, from
(\ref{definition-of-X-prime}) again, we have $X_{ij}^{\prime}=
-q^{-1}[i,1|j,n]X_{1n}^{-1}$ and
$X_{kl}^{\prime}=-q^{-1}[k,1|l,n]X_{1n}^{-1}$.  The desired relations
between $X_{ij}^{\prime}$ and $X_{kl}^{\prime}$ are deduced in an
obvious manner 
from \cite[Theorem 5.2.1]{PW} 
applied to the submatrix of $\mat{X}$ defined by 
rows $1,i,k$ and columns $j,l,n$.  \QED

The following proposition is a list of relations between entries of
$\widetilde{\mat{X}}$.  Relations between two elements of the first row
or of the last column of $\widetilde{\mat{X}}$ are deduced from the fact
that $\mat{X}$ is a $q$-quantum matrix. 

\begin{subproposition} \label{relations-in-X-tilda} 
The following relations hold between entries of the matrix
$\widetilde{{\bf X}}$.\\
\noindent 
1) If $1\leq j \leq n-1$ and $2\leq i \leq m$ then:
\[
X_{1j}X_{in}-q^{2}X_{in}X_{1j}=q(q^{2}-1)X_{ij}^{\prime}X_{1n}.
\]

\noindent 
2.1) For $1 \le j \le n-1$ and $2 \le k \le m$, 
\begin{eqnarray*} 
X_{1j}X_{kl}^{\prime}-X_{kl}^{\prime}X_{1j}&=&
(q^{-1}-q)X_{1l}X_{kj}^{\prime},\quad  \mbox{for}\; 
1 \le l \le j-1,\\
X_{1j}X_{kj}^{\prime}&=&q^{-1}X_{kj}^{\prime}X_{1j}, \\
X_{1j}X_{kl}^{\prime}&=&X_{kl}^{\prime}X_{1j}, \quad  \mbox{for}\; 
j+1 \le l \le n-1. 
\end{eqnarray*}

\noindent 
2.2) For $2 \le i \le m$ and $1 \le l \le n-1$,
\begin{eqnarray*}
X_{in}X_{kl}^{\prime}&=&X_{kl}^{\prime}X_{in},\quad \mbox{ for } \; 
2 \le k \le i-1,\\
X_{in}X_{il}^{\prime}&=&qX_{il}^{\prime}X_{in}, \\
X_{in}X_{kl}^{\prime}-X_{kl}^{\prime}X_{in}&=&
(q-q^{-1})X_{kn}X_{il}^{\prime}, \quad \mbox{ for } \;
i+1 \le k \le m.
\end{eqnarray*}

\noindent 
3) If $1 \le k < l \le n$ and $1 \le i < j \le m$ then
\[
X_{1k}X_{1l}=qX_{1l}X_{1k}
\quad \mbox{and} \quad X_{in}X_{jn}=qX_{jn}X_{in}.
\]
\end{subproposition}

\proof Part 3) is obvious.  The other relations are all obtained in the
following way.  First, use the definition 
(\ref{definition-of-X-prime}) to translate
the desired relation into a new one that uses 
only the entries of $\mat{X}$, by 
multiplying by $X_{1n}$.  The relation so obtained involves entries
coming from a certain submatrix of $\mat{X}$ of size at most $3 \times
3$.  Use Relation $\ref{2X2-q-matrix}$ and 
Relations $\ref{3X3-q-matrix}$ to check that the relation holds. 
Since these are easy but tedious computations we omit the details.  \QED

The next thing we want to do is to show that $\mat{X}'$ is a generic
$q$-quantum matrix.  We use a computation with Gelfand-Kirillov
dimension to do this.  The following lemma is probably well-known, but
we  have not located the exact statement that we need; and so we
include a proof.

\begin{sublemma} \label{first-GK-lemma} 
Let $B$ be a $\k$-algebra. 
Suppose $A$ is a subalgebra of $B$ and $x$ an element of $B$ such that
$B$ is generated by $A$ and $x$ as an algebra.  Further, suppose there
exists a finite dimensional subspace $V$ of $A$ that generates $A$ as an
algebra and such that $xV \subseteq Vx + A$. Then $\GK (B) \le \GK (A)
+1$.  \end{sublemma}

\proof Without loss of generality, we may assume that $1\in V$.  We
denote by $W$ the subspace of $B$ spanned by $x$ and $V$.  Thus, $W$ is
a finite dimensional subspace of $B$ that generates $B$ as an algebra. 
Since $xV \subseteq Vx + A$ and $V$ is finite dimensional, there exists
$m\in {\Bbb N}^\ast$ such that $xV \subseteq Vx + V^m$.  Then, an easy
induction shows that, for $n \in {\Bbb N}$, $x V^n \subseteq V^n x +
V^{m+n}$.  \\
We claim that,
\begin{equation} \label{incl}
W^{n} \subseteq V^{nm} +V^{nm} x+...+V^{nm} x^{n},
\end{equation} for all $n \in {\mathbb N}$. 
The inclusion (\ref{incl}) is trivial for $n=0$ 
(recall the standard convention that 
$V^0=W^0=\k$). Now, assume (\ref{incl}) holds for $p \in {\Bbb N}$, that is 
$W^{p} \subseteq V^{pm} +V^{pm} x+...+V^{pm} x^{p}$. It follows that
\begin{eqnarray*}
xW^{p}   & \subseteq &  xV^{pm} +xV^{pm} x+...+xV^{pm} x^{p} \\
         & \subseteq &  V^{pm}x + V^{pm+m} + V^{pm}x^{2} + V^{pm+m}x  
                + ...+V^{pm} x^{p+1} + V^{pm+m} x^{p}.
\end{eqnarray*} 
So $xW^{p} \subseteq V^{m(p+1)} +...+V^{m(p+1)} x^{p+1}$.
On the other hand, $VW^{p} \subseteq V^{mp+1} +...+V^{mp+1} x^{p}$.
We finally obtain 
$W^{p+1} =(V+\k x) W^{p} \subseteq V^{m(p+1)} +...+V^{m(p+1)} x^{p+1}$. 
This establishes (\ref{incl}) by induction on $n$. 

From (\ref{incl}), it follows that
$\dim W^{n} \le (n+1)\dim V^{mn}$, for all $n \in {\mathbb N}$, 
and thus
\[
\overline{\lim}\; \log_{n}\; \dim W^{n} \le
\overline{\lim}\; \log_{n}\; \dim V^{mn}+1.
\] 
This establishes that $\GK B \le \GK A +1$. \QED

\begin{subproposition} \label{X-prime-is-generic-q-quantum}
The matrix $\mat{X}'$ is a generic $q$-quantum matrix.
\end{subproposition}

\proof In this proof, we denote by $A$ the subalgebra of
${\k}_q[\mat{X}]_{X_{1n}}$ generated by the entries of $\mat{X}'$.  By
Lemma \ref{X-prime-is-q-quantum}, there is a surjective morphism $\phi
\; : \; {\cal O}_q(M_{m-1,n-1}) \longrightarrow A$.  We want to prove
that $\ker \phi = (0)$.  Assume that this is not the case.  Then,
since ${\cal O}_q(M_{m-1,n-1})$ is a domain, $\ker \phi$ must contain a
regular element of ${\cal O}_q(M_{m-1,n-1})$ and thus we have $\GK A <
\GK {\cal O}_q(M_{m-1,n-1}) = (m-1)(n-1)$. 

   By Lemma \ref{X-prime-is-q-quantum}, we see that the subalgebra $B$
of ${\k}_q[\mat{X}]_{X_{1n}}$ generated by 
$A$ and $X_{1n}$ is an extension
of $A$ of the type investigated in Lemma \ref{first-GK-lemma}.  Thus, we
must have $\GK B < (m-1)(n-1) + 1$.  Moreover, since $X_{1n}$ is central
in $B$, \cite[Proposition 4.2]{KL} 
shows that $\GK B_{X_{1n}} < (m-1)(n-1) + 1$.  Now, it
is clear that ${\k}_q[\mat{X}]_{X_{1n}}$ can be obtained by successive
algebra extensions starting from $B_{X_{1n}}$ and adding (in this
order) $X_{11},\dots,X_{1,n-1}$ and then $X_{mn},\dots,X_{2n}$ (this is
because ${\k}_q[\mat{X}]_{X_{1n}}$ is generated by the entries of
$\widetilde{\mat{X}}$ and $X_{1n}^{-1}$).  Moreover, Proposition
\ref{relations-in-X-tilda} shows that at each step, the extension is of
the type investigated in Lemma \ref{first-GK-lemma}.  Thus,
$(m-1)+(n-1)$ applications of this Lemma show that we must have $\GK
{\k}_q[\mat{X}]_{X_{1n}} < (m-1)(n-1) + 1 + (m-1)+(n-1) = mn$. 

However, we know that $\GK {\k}_q[\mat{X}]_{X_{1n}} \ge \GK
{\k}_q[\mat{X}] = mn$.  This is a contradiction and thus we have proved
that $\ker \phi = (0)$.  \QED 

Proposition~\ref{X-prime-is-generic-q-quantum} 
states that the subalgebra of
${\k}_q[\mat{X}]_{X_{1n}}$ generated by $\mat{X}'$ is isomorphic to
${\cal O}_q(M_{m-1,n-1})$;  for this reason, we denote it by
${\k}_q[\mat{X}']$.\\

The following remark will be useful in what follows.

\begin{subremark} \label{ore-ext-and-free-alg} 
\rm We denote by ${\rm F}\langle X_1,\dots,X_p\rangle $ the free
$\k$-algebra on $p$ generators $X_1,\dots,X_p$.  Let $I$ be the ideal
generated by elements $f_1,\dots,f_s \in {\rm F}\langle
X_1,\dots,X_p\rangle $, and set $A:={\rm F}\langle X_1,\dots,X_p\rangle
/I$. 
Finally, let $\sigma$ be an automorphism of ${\rm F}\langle
X_1,\dots,X_p\rangle $ and $\delta$ be a left $\sigma$-derivation of
${\rm F}\langle X_1,\dots,X_p\rangle $ such that $\sigma(I) = I$
and $\delta(I) \subseteq I$.  We denote by $\overline\sigma$ the
automorphism of $A$ induced by $\sigma$ and by $\overline\delta$ the
left $\overline\sigma$-derivation of $A$ induced by $\delta$.  Then,
$A[x;\overline\sigma,\overline\delta]$ is isomorphic to the algebra
$F\langle X_1,\dots,X_n,X\rangle /J$ where $J$ is the ideal of $F\langle
X_1,\dots,X_n,X\rangle $ generated by $f_1,\dots,f_s$ and the $p$
elements $XX_i - \sigma(X_i)X- \delta(X_i)$ (here we
identify ${\rm F}\langle X_1,\dots,X_p\rangle $ and the subalgebra of
${\rm F}\langle X_1,\dots,X_p,X\rangle $ generated by $X_1,\dots,X_p$). 
\end{subremark}

We now proceed to show that ${\k}_q[\mat{X}]_{X_{1n}}$ can be described
as an iterated Ore extension of its subalgebra ${\k}_q[\mat{X}']$
(recall 
the notation we fixed immediately after
Proposition~\ref{X-prime-is-generic-q-quantum}.  
For this, we first construct an algebra ${\cal B}$ which is an iterated
Ore extension of ${\cal O}_q(M_{m-1,n-1})$.  We then show that ${\cal
B}$ is isomorphic to ${\k}_q[\mat{X}]_{X_{1n}}$.  \\

Let us start with the algebra ${\cal O}_q(M_{m-1,n-1})$.  For
convenience of notation, we denote its canonical generators 
by $Y_{ij}^\prime$ for $2 \le i \le
m$ and $1 \le j \le n-1$.  Thus, $\mat{Y}'=
(Y_{ij}^\prime)_{2 \le i \le m,1 \le j \le n-1}$ is a generic
$q$-quantum matrix and, following our previous notation, we have ${\cal
O}_q(M_{m-1,n-1})=\k_q[\mat{Y}']$.  Now, we consider the Laurent
polynomial extension ${\cal A} = \k_q[\mat{Y}'][Y_{1n}^{\pm 1}]$
obtained from $\k_q[\mat{Y}']$ adding a central indeterminate denoted
$Y_{1n}$ as well as its inverse. 

   By Remark \ref{ore-ext-and-free-alg}, we know how to describe 
${\cal A}$ as a quotient of a free algebra.  Then, it is easy to see that
we can define a (unique) left skew derivation
$(\sigma_{11},\delta_{11})$ of ${\cal A}$ such that $\delta_{11}=0$ and
such that, $\sigma_{11}(Y_{1n})=qY_{1n}$ and, for $2 \le k \le m$ and $2
\le l \le n-1$, $\sigma_{11}(Y_{k1}^\prime)=q^{-1}Y_{k1}^\prime$ and
$\sigma_{11}(Y_{kl}^\prime)=Y_{kl}^\prime$.  We put ${\cal A}_{11}={\cal
A}[Y_{11};\sigma_{11},\delta_{11}]$. 

   Now, for $2 \le p \le n-1$, we construct an algebra ${\cal A}_{1p}$
by induction on $p$.  For $1 \le p \le n-2$, by Remark
\ref{ore-ext-and-free-alg}, we know how to describe ${\cal A}_{1p}$ as a
quotient of a free algebra.  It is then easy to check that we can define
a (unique) left skew derivation $(\sigma_{1,p+1},\delta_{1,p+1})$ of
${\cal A}_{1p}$ such that, for $2 \le k \le m$~:
\[
\begin{array}{lllll}
\sigma_{1,p+1}(Y_{kl}^\prime)=Y_{kl}^\prime &  
\mbox{ and } & \delta_{1,p+1}(Y_{kl}^\prime) = 
(q^{-1}-q)Y_{1l}Y_{k,p+1}^\prime & \mbox{ for } & 1 \le l \le p,\cr
\sigma_{1,p+1}(Y_{k,p+1}^\prime)=q^{-1}Y_{k,p+1}^\prime & \mbox{ and } &
\delta_{1,p+1}(Y_{k,p+1}^\prime)=0 & & \cr   
\sigma_{1,p+1}(Y_{kl}^\prime)=Y_{kl}^\prime & 
\mbox{ and } & \delta_{1,p+1}(Y_{kl}^\prime) = 0 & 
\mbox{ for } & p+2 \le l \le n-1 \cr   
\sigma_{1,p+1}(Y_{1n})=qY_{1n} & \mbox{ and } & 
\delta_{1,p+1}(Y_{1n})=0 & & \cr
\sigma_{1,p+1}(Y_{1l})=q^{-1}Y_{1l} & \mbox{ and } & 
\delta_{1,p+1}(Y_{1l})=0 & \mbox{ for } & 1 \le l \le p. \cr
\end{array}
\]
Then, for $1 \le p \le n-2$, we put ${\cal A}_{1,p+1}={\cal
A}_{1,p}[Y_{1,p+1};\sigma_{1,p+1},\delta_{1,p+1}]$. \\

Thus, we now have a first sequence of Ore extensions: ${\cal A}$, ${\cal
A}_{11}$,..., ${\cal A}_{1,n-1}$.  To finish the construction of ${\cal
B}$, we need a second such sequence which we now define.  \\

By Remark \ref{ore-ext-and-free-alg}, we know how to descibe 
${\cal A}_{1,n-1}$
as a quotient of a free algebra.  Then, it is easy to
see that we can define a (unique) left skew derivation
$(\sigma_{mn},\delta_{mn})$ of ${\cal A}_{1,n-1}$ such that, for $1 \le
l \le n-1$:
\[
\begin{array}{lllll}
\sigma_{mn}(Y_{kl}^\prime)=Y_{kl}^\prime &  
\mbox{ and } & \delta_{mn}(Y_{kl}^\prime) = 0  &
\mbox{ for } & 2 \le k \le m-1,\cr
\sigma_{mn}(Y_{ml}^\prime)=qY_{ml}^\prime & \mbox{ and } &
\delta_{mn}(Y_{ml}^\prime)=0 & & \cr   
\sigma_{mn}(Y_{1n})=q^{-1}Y_{1n} & \mbox{ and } & 
\delta_{mn}(Y_{1n})=0 & & \cr
\sigma_{mn}(Y_{1l})=q^{-2}Y_{1l} & \mbox{ and } & 
\delta_{mn}(Y_{1l})=(q^{-1}-q)Y_{ml}^\prime Y_{1n}. & & \cr
\end{array}
\]

We put ${\cal A}_{mn}={\cal A}_{1,n-1}[Y_{mn};\sigma_{mn},\delta_{mn}]$. \\

Now, for $2 \le s \le m-1$, we construct an algebra ${\cal A}_{sn}$, by {\it
decreasing} induction on $s$. 

For $3 \le s \le m$, by Remark \ref{ore-ext-and-free-alg}, 
we know how to describe
${\cal A}_{sn}$ as a quotient of a free algebra. 
It is then easy to check that  we can
define a (unique) left skew derivation 
$(\sigma_{s-1, n},\delta_{s-1,n})$ of ${\cal
A}_{sn}$ such that, for $1 \le l \le n-1$~:
\[
\begin{array}{lllll}
\sigma_{s-1,n}(Y_{kl}^\prime)= Y_{kl}^\prime &  
\mbox{ and } & \delta_{s-1,n}(Y_{kl}^\prime) = 0 &
\mbox{ for } & 2 \le k \le s-2  \cr
\sigma_{s-1,n}(Y_{s-1,l}^\prime) = qY_{s-1,l}^\prime & \mbox{ and } &
\delta_{s-1,n}(Y_{s-1,l}^\prime)= 0 & & \cr   
\sigma_{s-1,n}(Y_{kl}^\prime)=Y_{kl}^\prime & 
\mbox{ and } & \delta_{s-1,n}(Y_{kl}^\prime) = (q-q^{-1})Y_{kn}Y_{s-1,l}^\prime & 
\mbox{ for } &  s \le k \le m \cr   
\sigma_{s-1,n}(Y_{1n})=q^{-1}Y_{1n}, & \mbox{ and } & 
\delta_{s-1,n}(Y_{1n})=0 & & \cr
\sigma_{s-1,n}(Y_{1l})=q^{-2}Y_{1l}, & \mbox{ and } & 
\delta_{s-1,n}(Y_{1l})=(q^{-1}-q)Y_{s-1,l}^\prime Y_{1n} & &  \cr
\sigma_{s-1,n}(Y_{kn})=qY_{kn}, & \mbox{ and } & 
\delta_{s-1,n}(Y_{kn})=0 & \mbox{ for } & s \le k \le m. \cr
\end{array}
\]

Finally, we get a sequence of Ore extensions: ${\cal A}$, ${\cal
A}_{11}$,..., ${\cal A}_{1,n-1}$,${\cal A}_{mn}$,...,${\cal A}_{2n}$,
and we put ${\cal B}= {\cal A}_{2n}$.  Iterative applications of Remark
\ref{ore-ext-and-free-alg} show that ${\cal B}$ can be easily described
as the quotient of a free algebra in $mn$ generators.  Clearly, we have
constructed ${\cal B}$ in such a way that the relations between the $mn$
generators of this algebra are exactly the same as those holding
between the generators of $\widetilde{\mat{X}}$ 
as noted in Proposition 1.1.2.  It follows at once
that we can define a morphism of algebras 
\[ 
\begin{array}{ccrclllll}
\varphi & : & {\cal B} & \longrightarrow & \k_q[\mat{X}]_{X_{1n}} \cr
 & & Y_{kl}^\prime & \mapsto & X_{kl}^\prime & \mbox{ for } 
& 2 \le k \le m & \mbox{
and } & 1 \le l \le n-1 \cr
 & & Y_{1n} & \mapsto & X_{1n} & & & & \cr
 & & Y_{1l} & \mapsto & X_{1l} & \mbox{ for } & 1 \le l \le n-1 & & \cr
 & & Y_{kn} & \mapsto & X_{kn} & \mbox{ for } & 2 \le k \le m & &  \cr
\end{array}
\]

\begin{subproposition} \label{ore-ext-iso} 
The morphism $\varphi$ is an isomorphism.
\end{subproposition}

\proof The surjectivity of $\varphi$ is obvious since the entries of
$\widetilde{\mat{X}}$ together with $X_{1n}^{-1}$ form a set of algebra
generators for $\k_q[\mat{X}]_{X_{1n}}$.  It remains to prove that $\ker
\varphi = (0)$.  Recall that ${\cal B}$ is obtained from ${\cal A}={\cal
O}_q(M_{m-1,n-1})[X_{1n}^{-1}]$ by $(m-1) + (n-1)$ succesive Ore extensions
that are all extensions of algebras of the type investigated in Lemma
\ref{first-GK-lemma}.  Thus, we have $\GK {\cal B} \le \GK {\cal A} +
(m-1) + (n-1)$.  On the other hand, ${\cal A}$ is a Laurent polynomial
extension of ${\cal O}_q(M_{m-1,n-1})$ thus $\GK {\cal A} = \GK {\cal
O}_q(M_{m-1,n-1}) +1 = (m-1)(n-1) + 1$, by \cite[Proposition 3.5 and
Proposition 4.2]{KL}.  All this 
together gives
\[
\GK {\cal B} \le (m-1)(n-1) + 1 + (m-1) + (n-1) = mn.
\]
On the other hand, 
we have $\GK \k_q[\mat{X}]_{X_{1n}} \ge \GK
\k_q[\mat{X}] = mn$, since $\k_q[\mat{X}]$ is a subalgebra of
$\k_q[\mat{X}]_{X_{1n}}$. 

 Now, ${\cal B}$ is clearly a domain; and so if we assume that $\ker
\varphi \neq 0$ then $\GK
{\cal B}/\ker\varphi < \GK {\cal B} = mn$.  But this is a contradiction
since ${\cal B}/\ker\varphi \cong \k_q[\mat{X}]_{X_{1n}}$.  Thus,
$\varphi$ is injective.  \QED

\section{Reduction of the size of quantum minors}

Using the previous section we are now able to link $k \times k$
quantum minors 
of ${\bf X}$, for $k\geq 2$, 
with $(k-1) \times (k-1)$ quantum minors of ${\bf
X}^{\prime}$ (provided $X_{1n}$ is invertible, of course). \\

Recall from the introduction that the expression $[i_{k}, \dots
,i_{1}|j_{1}, \dots ,j_{k}]$ stands for a $k \times k$ quantum minor of
$\mat{X}$, namely the quantum determinant of the submatrix of ${\bf X}$
obtained from $\mat{X}$ using rows $i_{1},\dots ,i_{k}$ and
columns $j_{1},\dots ,j_{k}$.  This notation is extended to ${\bf
X}^{\prime}$ adding a ``$\,^{\prime}\,$'' to avoid confusion. Thus, a
quantum $k \times k$ minor of ${\bf X}^{\prime}$ will be expressed by a
symbol $[i_{k},\dots ,i_{1}|j_{1},\dots ,j_{k}]^{\prime}$.  The
convention on the index of rows of ${\bf X}^{\prime}$ (see introduction
of Section~\ref{newsetofgen}) is in order; so,   
in such an expression, we shall
always have $i_{1} \ge 2$ and $j_{k} \le n-1$.  \\

The results we need will follow from the special case where $m=n$.  In
this context, the role played by $(n-1) \times (n-1)$ minors is of
special importance; thus we use a more convenient notation for them
(coming from \cite{PW}).  For $1\le i,j \le n$, the $(n-1) \times (n-1)$
quantum minor of $\mat{X}$ obtained by deleting the i-th row and j-th
column is denoted $A(ij)$.  Moreover, for $2 \le i \le n$ and $1 \le j
\le n-1$ we also define $A'(ij)$ to be the $(n-2) \times (n-2)$ quantum
minor of $\mat{X}^{\prime}$ obtained from $\mat{X}^{\prime}$ by deleting 
the i-th row and j-th column.\\

\begin{subtheorem} \label{reduction-theorem} Assume that  $m=n$.  
Then, with the
above notation: 
\[
(\det_{q}{\bf X}^{\prime}) X_{1n}=X_{1n} (\det_{q}{\bf
X}^{\prime}) =(-q)^{1-n}\det_{q}{\bf X}.
\]
  \end{subtheorem}

\proof Note that $(\det_{q}{\bf X}^{\prime}) X_{1n}=X_{1n} (\det_{q}{\bf
X}^{\prime})$ is clear from \ref{X-prime-is-q-quantum}.  The proof is by
induction on $n$.  The case where $n=2$ is an obvious consequence of
(\ref{definition-of-X-prime}).  We suppose now that the result is true
for any integer less than or equal to $n-1$.  Because of Lemma 
\ref{X-prime-is-q-quantum} the relations of \cite[Corollary 4.4.4]{PW} 
give us the
expansion: 
\begin{equation} {\rm det} _{q}{\bf X}^{\prime} =
\sum\limits_{j=1}^{n-1} (-q)^{j-1}X_{2j}^{\prime} A'(2j). 
\end{equation} 
Now, 
$X_{1n}X_{2j}^{\prime} =
X_{2j}^{\prime}X_{1n}=X_{2j}X_{1n}-q^{-1}X_{1j}X_{2n}$ for all $j \in
\{1,\dots,n-1\}$; and so 
\[
\begin{array}{rcl} X_{1n}\det _{q}{\bf X}^{\prime}  &=& 
\sum\limits_{j=1}^{n-1} (-q)^{j-1}X_{1n}X_{2j}^{\prime} A'(2j) \\ 
&=& \sum\limits_{j=1}^{n-1} 
(-q)^{j-1}(X_{2j}X_{1n}-q^{-1}X_{1j}X_{2n})A'(2j) \\ 
&=& \sum\limits_{j=1}^{n-1} (-q)^{j-1}X_{2j}X_{1n} A'(2j) 
-q^{-1}\sum\limits_{j=1}^{n-1} (-q)^{j-1}X_{1j}X_{2n} A'(2j). 
\end{array}
\] 
By the induction hypotheses we have 
$X_{1n}A'(2j) = (-q)^{2-n}A(2j)$  
for all $j \in \{1,\dots,n-1\}$. Hence, 
\[
\begin{array}{rcl}
X_{1n}{\rm det}_{q}{\bf X}^{\prime} &=&
\sum\limits_{j=1}^{n-1}(-q)^{j-1}X_{2j}(-q)^{2-n}A(2j)
-q^{-1}\sum\limits_{j=1}^{n-1}
(-q)^{j-1}X_{1j}X_{2n}X_{1n}^{-1}(-q)^{2-n}A(2j)\\
&=&
\sum\limits_{j=1}^{n-1}(-q)^{j+1-n}X_{2j}A(2j)
-q^{-1}\sum\limits_{j=1}^{n-1}(-q)^{j+1-n}X_{1j}X_{2n}X_{1n}^{-1}A(2j)
\end{array}
\] 
Setting 
$R=-q^{-1}\sum\limits_{j=1}^{n-1}
(-q)^{j+1-n}X_{1j}X_{2n}X_{1n}^{-1}A(2j)$, we
obtain:\\
\[
\begin{array}{rcl} X_{1n}R &=& 
\sum\limits_{j=1}^{n-1}(-q)^{j-n}X_{1j}X_{2n}A(2j) \\ &=& 
\sum\limits_{j=1}^{n-1}(-q)^{j-n}
(X_{2n}X_{1j}+(q-q^{-1})X_{1n}X_{2j})A(2j) \\ &=&
\sum\limits_{j=1}^{n-1}(-q)^{j-n}X_{2n}X_{1j}A(2j) + (q-q^{-1})
\sum\limits_{j=1}^{n-1}(-q)^{j-n}X_{1n}X_{2j}A(2j). \\
\end{array}
\] 
On the other hand, from \cite[Corollary 4.4.4]{PW}, 
we have the relation
$\sum\limits_{j=1}^{n} (-q)^{j-2}X_{1j}A(2j)=0$. Thus,
\[
\begin{array}{rcl} X_{1n}R &=& 
(-q)^{2-n}X_{2n}\sum\limits_{j=1}^{n-1}(-q)^{j-2}X_{1j}A(2j) + (q-q^{-1})
(-q)^{2-n}X_{1n}\sum\limits_{j=1}^{n-1}(-q)^{j-2}X_{2j}A(2j). \\ &=& 
(-q)^{2-n}X_{2n}(0-(-q)^{n-2}X_{1n}A(2n)) + (q-q^{-1})
(-q)^{2-n}X_{1n}\sum\limits_{j=1}^{n-1}(-q)^{j-2}X_{2j}A(2j). \\
\end{array}
\] 
It follows that: 
\[
R = (-q)^{-1}X_{2n}A(2n) +
(q-q^{-1})(-q)^{2-n}\sum\limits_{j=1}^{n-1}(-q)^{j-2}X_{2j}A(2j).
\]
  Thus, 
\[
\begin{array}{rcl} X_{1n} \det _{q} {\bf X}^{\prime} 
&=& \sum\limits_{j=1}^{n-1}
(-q)^{j+1-n}X_{2j}A(2j)  +(-q)^{-1}X_{2n}A(2n) \\ &&
\;\;\;\;\;\;\;\;\;\;\;\;\;\;\;\;\;\;\;+ (q-q^{-1})
(-q)^{2-n}\sum\limits_{j=1}^{n-1}(-q)^{j-2}X_{2j}A(2j)\\ 
&=&
(-q)^{2-n}\sum\limits_{j=1}^{n-1} ((-q)^{j-1}+(q-q^{-1})(-q)^{j-2})X_{2j}A(2j)
+(-q)^{-1}X_{2n}A(2n) \\ 
&=& (-q)^{2-n}\sum\limits_{j=1}^{n-1}
((-q)^{j-1}-(-q)^{j-1}-q^{-1}(-q)^{j-2})X_{2j}A(2j) +(-q)^{-1}X_{2n}A(2n) \\ 
&=& (-q)^{1-n}\sum\limits_{j=1}^{n-1} (-q)^{j-2}X_{2j}A(2j)
+(-q)^{1-n}(-q)^{n-2}X_{2n}A(2n) \\ 
&=&
(-q)^{1-n}\sum\limits_{j=1}^{n} (-q)^{j-2}X_{2j}A(2j).
\end{array}
\] 
Again, by using \cite[Corollary 4.4.4]{PW}, it follows that  
$X_{1n} \det _{q} {\bf X}^{\prime} = (-q)^{1-n}\det _{q} {\bf X}$. \QED

By using Theorem \ref{reduction-theorem}, we can establish Corollary
\ref{cor-reduction-theorem} which links $(p-1) \times (p-1)$ minors of
$\mat{X}^\prime$ with $p \times p$ minors of $\mat{X}$, for $p\geq2$, 
that involve the
first row and the last column of $\mat{X}$. 

\begin{subcorollary} \label{cor-reduction-theorem} Let $p \ge 2$.
Suppose that 
$I=\{1= i_{1}<\dots<i_{p}\} \subseteq {\Bbb N}_m$ and
$J=\{j_{1}<\dots<j_{p}=n\} \subseteq {\Bbb N}_n$ and 
set  $I'=\{i_{2}<\dots<i_{p}\}$ and
$J'=\{j_{1}<\dots<j_{p-1}\}$. Then
\[
[I'|J']'=(-q)^{1-p}[I|J]X_{1n}^{-1}=(-q)^{1-p}X_{1n}^{-1}[I|J].  
\]
\end{subcorollary}

\proof This is an immediate consequence of 
Theorem \ref{reduction-theorem} applied
to the square submatrix obtained from $\mat{X}$ by using rows
$i_{1},\dots,i_{p}$ and columns $j_{1},\dots,j_{p}$ of $\mat{X}$.  
\qed\\

Recall from the introduction that, for $1 \le t \le \min\{m,n\}$, we
denote by ${\cal I}_{t}(\mat{X})$ the ideal generated in
$\k_{q}[\mat{X}]$ by the $t \times t$ minors of $\mat{X}$.  Clearly, the
ideal generated by the $t \times t$ minors of $\mat{X}$ in
$\k_{q}[\mat{X}]_{X_{1n}}$ is just $\widetilde{{\cal
I}}_{t}(\mat{X}):={\cal I}_{t}(\mat{X})\k_{q}[\mat{X}]_{X_{1n}}$. 

\begin{sublemma} \label{smaller-set-of-gen} Let $1 \le t \le
\min\{m,n\}$.  \\ (i) ${\cal I}_{t}(\mat{X})$ coincides with the right
ideal of $\k_{q}[\mat{X}]$ generated by the $t \times t$ minors of
$\mat{X}$.  \\ (ii) $\widetilde{{\cal I}}_{t}(\mat{X})$ coincides with
the right ideal of $\k_{q}[\mat{X}]_{X_{1n}}$ generated by the $t \times
t$ minors $[I|J]$ of $\mat{X}$ such that $1 \in I$ and $n \in J$.  \\
\end{sublemma}

\proof (i) The case where $m=n$ follows at once from \cite[Corollary
A.2]{GL}.  
Now, set $s:=\max\{m,n\}$.  There is a surjective algebra
morphism
\[
\begin{array}{ccrcl}
\pi & : & {\cal O}_{q}(M_s) & \longrightarrow & {\cal O}_q(M_{m,n}) \cr
 & & X_{ij} & \mapsto & \left\{
              \begin{array}{ll} X_{ij} & \mbox{ if } i \le m 
\mbox{ and } j \le n \cr
                                    0  & \mbox{ otherwise }
              \end{array} \right.
\end{array}
\]
and, for $I,J \subseteq {\Bbb N}_s$ such that $|I|=|J|=t$, a $t \times
t$ minor $[I|J]$ of ${\cal O}_q(M_s)$ is sent to $[I|J]$ if $I \subseteq
{\Bbb N}_m$ and $J \subseteq {\Bbb N}_n$ and is sent to $0$ otherwise. 
It follows that the ideal of ${\cal O}_q(M_{m,n})$ generated by the $t
\times t$ minors is the image under $\pi$ of the ideal of ${\cal
O}_q(M_{s})$ generated by the $t \times t$ minors.  From this, we see
that point (i) for arbitrary positive integers $m$ and $n$ follows from
the special case $m=n$.  

(ii) By part (i), $\widetilde{{\cal
I}}_{t}(\mat{X})$ coincides with the right ideal of
$\k_{q}[\mat{X}]_{X_{1n}}$ generated by the $t \times t$ minors.  For
the purpose of this proof, denote by ${\cal S}$ the set of $t
\times t$ minors $[I|J]$ of $\mat{X}$ such that $1 \in I$ and $n \in J$. 

Let $[I|J]$ be a $t \times t$ minor of $\mat{X}$ such that $1 \in I$ but
$n \not\in J$.  We may apply \cite[Corollary 4.4.4]{PW} 
to the subalgebra of
$\k_q[\mat{X}]$ generated by those $X_{ij}$ such that $i \in I$ and $j
\in J \cup \{n\}$.  This leads to the equation:
\[
\sum_{j \in J\cup\{n\}} (-q)^\bullet [I|J_j]X_{1j}=0,
\]
where, for $ j \in J \cup \{n\}$, we put $J_j = J \cup \{n\} \setminus
\{j\}$, and occurences of $(-q)^\bullet$ denote integer powers of $-q$
which it is not necessary to specify exactly.  
Since $J_n=J$, it follows that the equation
\[
[I|J]=-\sum_{j \in J} (-q)^\bullet [I|J_j]X_{1j}X_{1n}^{-1}
\]
holds in $\k_q[\mat{X}]_{X_{1n}}$.  
Now, $n \in J_j$ for each $j \in J$; and so 
we have shown that $[I|J]$ is in the right ideal generated by ${\cal S}$
in $\k_q[\mat{X}]_{X_{1n}}$.  By a similar argument, we prove that a $t
\times t$ minor $[I|J]$ such that $1 \not\in I$ but $n \in J$ is in the
right ideal generated by ${\cal S}$ in $\k_q[\mat{X}]_{X_{1n}}$. 

It remains to deal with a $t \times t$ minor $[I|J]$ such that $1 \notin
I$ and $n \notin J$.  In this case, we apply \cite[Corollary 4.4.4]{PW}
in 
the subalgebra of $\k_q[\mat{X}]$ generated by those $X_{ij}$ such that
$i \in I \cup \{1\}$ and $j \in J \cup \{n\}$.  This gives us the
relation 
\begin{equation} \label{eq-1} 
[I \cup \{1\}|J \cup
\{n\}]=\sum_{j \in J \cup \{n\}} (-q)^\bullet [I|J_j]X_{1j} =\sum_{j \in
J} (-q)^\bullet [I|J_j]X_{1j} + (-q)^\bullet [I|J]X_{1n} 
\end{equation}
where, for $ j \in J \cup \{n\}$, we put $J_j = J \cup \{n\} \setminus
\{j\}$.  
Another application of \cite[Corollary 4.4.4]{PW} gives
\begin{equation} \label{eq-2} 
[I \cup \{1\}|J \cup \{n\}]=\sum_{j \in J
\cup \{n\}} (-q)^\bullet [I\cup\{1\}\setminus \{s\}|J_j]X_{sj},
\end{equation} 
where, $s=\max I$.  Equation (\ref{eq-2}) and the results
we established above show that the $(t+1) \times (t+1)$ minor $[I \cup
\{1\}|J \cup \{n\}]$ is in the right ideal generated in
$\k_q[\mat{X}]_{X_{1n}}$ 
by ${\cal S}$.  Thus, using (\ref{eq-1}), it follows that $[I|J]$ is
also in the right ideal generated in $\k_q[\mat{X}]_{X_{1n}}$ 
by ${\cal S}$.  The
proof is now complete.  \qed\\

Recall that, by Lemma \ref{X-prime-is-q-quantum} and Proposition 
\ref{X-prime-is-generic-q-quantum} we know that the subalgebra of
$\k_q[\mat{X}]_{X_{1n}}$ generated by the $X_{ij}^\prime$ for $2 \le i
\le m$ and $1 \le j \le n-1$ is isomorphic to ${\cal O}_q(M_{m-1,n-1})$
and that we denote it by $\k_q[\mat{X}^\prime]$.  Following our
conventions, if $2 \le t \le \min\{m,n\}$, we denote by ${\cal
I}_{t-1}(\mat{X}^\prime)$ the ideal of $\k_q[\mat{X}^\prime]$ generated
by the $(t-1) \times (t-1)$ minors of $\mat{X}^\prime$.  In this
notation, we have the following important result. 

\begin{subproposition} \label{main-consequence-of-reduction} 
For $2 \le t \le
\min\{m,n\}$, the following equality holds: 
\[
\widetilde{{\cal I}}_{t}(\mat{X})=
{\cal I}_{t-1}(\mat{X}^{\prime})\k_{q}[\mat{X}]_{X_{1n}}.
\]
\end{subproposition}

\proof By Lemma \ref{smaller-set-of-gen} (ii), 
we know that $\widetilde{{\cal
I}}_{t}(\mat{X})$ coincides with the right ideal of
$\k_{q}[\mat{X}]_{X_{1n}}$ generated by the $t \times t$ minors $[I|J]$
of $\mat{X}$ such that $1 \in I$ and $n \in J$.  On the other hand, let
$[I|J]$ be a $t \times t$ minor of $\mat{X}$ such that $1 \in I$ and $n
\in J$.  
Applying Corollary~\ref{cor-reduction-theorem}, we have 
$(-q)^{1-t}[I|J]=[I\setminus\{1\}|J\setminus\{n\}]'X_{1n}$.  Thus,
$\widetilde{{\cal I}}_{t}(\mat{X})$ coincides with the right ideal of
$\k_{q}[\mat{X}]_{X_{1n}}$ generated by the $(t-1) \times (t-1)$ minors
$[I|J]'$ of $\mat{X}'$ (such that $I \subseteq \{2,\dots,n\}$ and $J
\subseteq \{1,\dots,m-1\}$).  On the other hand, Lemma 
\ref{smaller-set-of-gen} (i) shows that ${\cal
I}_{t-1}(\mat{X}^{\prime})$ is the right ideal of $\k_{q}[\mat{X}']$
generated by the $(t-1) \times (t-1)$ minors $\mat{X}'$; so, the proof
is complete.  \qed\\

  Recall from the introduction that, for any positive integer $t$ such
that $t \le \min \{m,n\}$, we define the quantum determinantal ring
$R_{t}(\mat{X}):=\k_{q}[\mat{X}]/{\cal I}_{t}(\mat{X})$. This is a domain
by \cite[Corollary 2.6]{GL}.  
If we put $x_{ij}:=X_{ij} + {\cal I}_{t}(\mat{X})$, 
for $1 \le i \le m$ and $1 \le j \le n$, then there is a canonical
isomorphism 
\begin{equation} R_{t}(\mat{X})_{x_{1n}} \cong
\k_{q}[\mat{X}]_{X_{1n}}/\widetilde{{\cal I}}_{t}(\mat{X}). 
\end{equation} 
We finish this subsection by showing that, for $2 \le t
\le \min\{m,n\}$, the ring 
$R_{t}(\mat{X})_{x_{1n}}$ can be described as an iterated Ore
extension of $R_{t-1}(\mat{X}')$.  To achieve this aim, we will have to
make use of Proposition \ref{ore-ext-iso} which shows that
$\k_q[\mat{X}]_{X_{1n}}$ can be obtained from its subalgebra
$\k_q[\mat{X}^\prime]$ by iterated Ore extensions adding successively
$X_{1n}^{\pm 1}$, $X_{11}$,...,$X_{1,n-1}$,$X_{mn}$,...,$X_{2n}$. 

\begin{subtheorem} \label{iterated-ore-ext} 
For $2 \le t \le \min\{m,n\}$, the ring 
$R_{t}(\mat{X})_{x_{1n}}$ is a localisation of an iterated 
Ore extension of $R_{t-1}(\mat{X}')$.
\end{subtheorem}

\proof We start by giving a list of relations between a $(t-1) \times
(t-1)$ minor of $\mat{X}'$ and the generators $X_{1n}^{\pm 1}$,
$X_{11}$,...,$X_{1,n-1}$,$X_{mn}$,...,$X_{2n}$.  (Note that, by Lemma
\ref{X-prime-is-q-quantum}, we already know that a $(t-1) \times (t-1)$
minor of $\mat{X}'$ commutes with $X_{1n}^{\pm 1}$.)

   Let $I' \subseteq \{2,\dots,m\}$ and $J' \subseteq \{1,\dots,n-1\}$
be sets of indices such that $|I'|=|J'|=t-1$. 
Setting $I=I' \cup \{1\}$
and $J=J' \cup \{n\}$, Corollary~\ref{cor-reduction-theorem}, shows that
$(-q)^{1-t}[I|J]=[I'|J']'X_{1n}$.  
On the other hand, \cite[Lemma 4.5.1 and Theorem 
4.6.1]{PW} 
give the following relations, for $k \in \{2,\dots,m\}$ and $l
\in \{1,\dots,n-1\}$:

\Q 
1) $X_{1l}[I|J]=[I|J]X_{1l}$ if $l \in J'$; \\
\hspace{5ex}
2) $X_{1l}[I|J]-q[I|J]X_{1l} =q(q-q^{-1})\sum_{j<l,j\in J'}(-q)^\bullet
X_{1j}[I|J\cup\{l\}\setminus\{j\}]$ if $l \not\in J'$; \\ 
\hspace{5ex}
3) $X_{kn}[I|J]=[I|J]X_{kn}$ if $k \in I'$; \\ 
\hspace{5ex} 
4) $X_{kn}[I|J]-q^{-1}[I|J]X_{kn}=q^{-1}(q^{-1}-q)\sum_{j>k,j\in
I'}(-q)^\bullet X_{jn}[I\cup\{k\}\setminus\{j\}|J]$ if $k \not\in I'$.

\Q 
(Points 1) and 3) follow from \cite[Theorem 4.6.1]{PW}, point 2) is 
\cite[Lemma 4.5.1(1) first relation]{PW},  
point 4) is \cite[Lemma 4.5.1(3) second relation]{PW}.) Hence, for $k
\in \{2,\dots,n\}$ and $l \in \{1,\dots,n-1\}$, 
Corollary~\ref{cor-reduction-theorem} gives: 

\Q 
1')
$X_{1l}[I'|J']'=q^{-1}[I'|J']'X_{1l}$ if $l \in J'$; \\
 2')
$X_{1l}[I'|J']'-[I'|J']'X_{1l} =q(q-q^{-1})\sum_{j<l,j\in
J'}(-q)^\bullet X_{1j}[I'|J'\cup\{l\}\setminus\{j\}]'$ if $l \not\in J'$;
\\
 3') $X_{kn}[I'|J']'=q[I'|J']'X_{kn}$ if $k \in I'$; \\
 4')
$X_{kn}[I'|J']'-[I'|J']'X_{kn}=q^{-1}(q^{-1}-q)\sum_{j>k,j\in
I'}(-q)^\bullet X_{jn}[I'\cup\{k\}\setminus\{j\}|J']'$ if $k \not\in
I'$.  

\Q 
 By Proposition \ref{ore-ext-iso}, $\k_q[\mat{X}]_{X_{1n}}$ can
be obtained from its subalgebra $\k_q[\mat{X}^\prime]$ by iterated Ore
extensions adding successively $X_{1n}^{\pm 1}$,
$X_{11}$,...,$X_{1,n-1}$,$X_{mn}$,...,$X_{2n}$.  The relations 1') to
4') together with the fact that $X_{1n}$ commutes with any element of
the subalgebra $\k_q[\mat{X}']$ show that, at each step of this Ore
extension, the ideal of the base algebra generated by the $(t-1) \times
(t-1)$ minors of $\mat{X}'$ is invariant under the corresponding skew
derivation.  It follows that $\k_q[\mat{X}]_{X_{1n}}/{\cal
I}_{t-1}(\mat{X}')\k_q[\mat{X}]_{X_{1n}}$ is isomorphic to an iterated
Ore extension of $\k_q[\mat{X}']/{\cal I}_{t-1}(\mat{X}')$.  But
Proposition \ref{main-consequence-of-reduction} shows that
$\widetilde{{\cal I}}_{t}(\mat{X})= {\cal
I}_{t-1}(\mat{X}^{\prime})\k_{q}[\mat{X}]_{X_{1n}}$.  The proof is thus
complete.  \QED

\section{Quantum determinantal rings are maximal orders}

Let $R$ be a noetherian domain with division ring of fractions $Q$. 
Then $R$ is said to be a {\bf maximal order} in $Q$ if the following
condition is satisfied: if $T$ is a ring such that $R\subseteq T
\subseteq Q$ and such that there exist nonzero elements $a, b \in R$
with $aTb\subseteq  R$, then $T=R$.  This condition is the natural
noncommutative analogue of normality for commutative domains, see, for
example, \cite[Section 5.1]{McCR}.  In this section, we investigate the
maximal order condition for determinantal rings.

Recall that the quantum determinantal ring
$R_{t}(\mat{X}):=\k_{q}[\mat{X}]/{\cal I}_{t}(\mat{X})$ is a domain by
\cite[Corollary 2.6]{GL}.  Also, recall that $R_2(\mat{X})$ has been
shown to be a maximal order in \cite{R}. 
In this section, we prove that $R_{t}(\mat{X})$ is
a maximal order in its division ring of fractions,  
when ${\mathbb K}={\mathbb C}$ and for 
$q$ an element of ${\mathbb C}$ transcendental over ${\mathbb Q}$.

The following Lemma from \cite{R} (see \cite[Lemma 1.1]{R}) 
is recalled here for the convenience of the reader. 

\begin{lemma} \label{lemme-mo} 
Let $R$ be a noetherian domain and
$Q={\rm Frac}\,R$ its division ring of fractions.  Assume there exists a
nonzero normal element $x$ in $R$ such that 

1) $xR=\cap_{i=1}^{r} {\goth
p_{i}}$ where, for $i \in \{1,\dots,r\}$, ${\goth p}_{i}$ is a
completely prime ideal of $R$,  and 

2) the localisation $R_{x}$ of $R$
with respect to the set $\{x^i,\,i \in {\Bbb N}\}$ is a maximal order in
$Q$.   

\noindent 
Let $\tau$ be the automorphism of $R$ associated with $x$; that is 
$ax = x\tau(a)$ for all $a\in R$.  
Suppose that $\tau({\goth p}_{i}) \subseteq {\goth p}_{i}$ for $i
\in \{1,\dots,r\}$.  Then $R$ is a maximal order in  $Q$.  
\end{lemma}

One case where the above result immediately applies is the case that the
ideal $xR$ itself is a completely prime ideal. It is this case that we
want to use. The applicability of the above lemma to determinantal rings
is a consequence of the results of the previous section, since we can
deduce the following result.

\begin{theorem}\label{main-thm}
Assume that $\k = {\mathbb C}$ and let 
$q$ be an element of ${\mathbb C}$ transcendental over ${\mathbb
Q}$. Let $t$ be an integer
such that $0 < t \leq \min\{m,n\}$.  Then $R_{t}(\mat{X})$ is a maximal
order. 
\end{theorem} 

\proof For $t =1$, this is trivial and for $t=2$ it is
\cite[Th\'eor\`eme 2.3.11]{R}. We proceed by induction on $t$. Assume
that the result is true for an integer $s \geq 2$ and set $t= s+1$. By
Theorem~\ref{iterated-ore-ext}, $R_{t}(\mat{X})_{x_{1n}}$ is a
localisation of an iterated Ore extension of $R_{t-1}(\mat{X}')$. The
induction hypothesis shows that $R_{t-1}(\mat{X}')$ is a maximal order
in its quotient ring and so $R_{t}(\mat{X})_{x_{1n}}$ is a maximal order
by \cite[V.2.5 and IV.2.1]{MR}. However, \cite[Corollaire 11.7]{Lau}
shows that $\langle x_{1n}\rangle$ is a completely prime ideal of
$R_{t}(\mat{X})$. Thus, Lemma~\ref{lemme-mo} shows that $R_{t}(\mat{X})$
is a maximal order. \QED


We conjecture that this result holds for arbitrary non-zero $q$ in any
field $\k$ and for all $t$.  We conclude by proving that the factor by
the quantum determinantal ideal $\langle \det_q \rangle$ in $\oqmn$ is a
maximal order.  All that remains to be proved after the above discussion
is that the ideal $\langle \det_q , X_{1n}\rangle$ is a completely prime
ideal.  This is what we do next.

\section{$\oqmn/\langle \det_q \rangle$ is a maximal order} 

In this section  we need
to use the {\bf preferred bases} in ${\cal O}_q(M_{u,v})$ 
developed in \cite{GL}, and we follow
the notation of that paper. See, in particular, \cite[Corollary
1.11]{GL}. We recall the
notation $[T|T']$ for the product of quantum minors corresponding to an
allowable bitableau $(T,T')$.  We recall also that it is sometimes
convenient to label rows of $(T,T')$ in the form $(I,J)$ where $I$ and
$J$ are sets of row and column indices, respectively (of course, $I
\subseteq \{1,\dots,u\}$ and $J \subseteq \{1,\dots,v\}$); such a pair 
is called an {\bf index pair} (see \cite[Section 1]{GL}).  Many of
the results in \cite{GL} are stated for the square case ${\cal
O}_q(M_{u,u})$, and there are easy extensions to the rectangular case   
${\cal O}_q(M_{u,v})$, see, for example, \cite[1.11]{GL} and
\cite[Section 2]{GLR} for more details of this standard procedure.

We aim to prove that the ideal $\langle \det_q , X_{1n}\rangle$ of
$\oqmn$ is a completely prime ideal for each $n\geq 3$, by using the
following result of David Jordan, \cite{J}. 

\begin{proposition} 
Let $\sigma$ be an automorphism  and $\delta$ be a $\sigma$-derivation
of a domain $A$. Let $R=A[x;\sigma,\delta]$. Let $c$
be a normal element of $R$ of the form $dx+e$, where $d,e\in A$
and $d\neq 0$.
Let $\beta$ be the automorphism of
$R$ such that $cr=\beta(r)c$ for all $r\in R$. 
Then $\beta(A)=A$, the element 
$d$ is normal in $A$ and $\beta(a)d=d\sigma(a)$ for
all
$a\in A$. Furthermore, if $e$ is regular modulo the ideal $Ad = dA$
then
$R/Rc$ is a domain.
\label{Jordan}
\end{proposition}

\begin{proposition} \label{prop-dom} 
The ideal $\langle \det_q , X_{1n}\rangle$ of
$\oqmn$ is a completely prime ideal for each $n\geq 3$.  
\end{proposition}

\proof 
The idea is first to factor out $X_{1n}$ obtaining a domain and then to
factor out $D_n := \det_q$ and see that we still have a 
domain by invoking Jordan's
Domain Theorem.  We use the notation from that theorem.  Let $T$ be the
subalgebra of $\oqmn$ generated by all $X_{ij}$ except $X_{nn}$.  Set
$A:= T/\ideal{X_{1n}}$.  Then $A$ is a domain.  Let $R:= A[X_{nn};
\sigma, \delta] \cong
{\cal O}_q(M_n)/\ideal{X_{1n}}$.  
Let $D_n$ denote the quantum determinant of
$\oqmn$ and let $D_{n-1}$ denote the quantum determinant of the copy of
${\cal O}_q(M_{n-1})$ generated by the $X_{ij}$ with $i, j < n$.  Note that
$D_{n-1} = A(nn)$ in the notation introduced before
Theorem~\ref{reduction-theorem}.  Set $x:= X_{nn}
\in R$ and $d:= \widebar{D_{n-1}} \in A \subseteq R$.  Note that $d \neq
0$ in $R$.  Finally, set $c := \widebar{D_n} \in R$.  The quantum
Laplace expansion of $D_n$ by the $n$th column gives 
$D_n = \sum_{i=1}^{n}\, \pm\qdot A(in)X_{in}$ so that $c = dx + e$,
where $e$ is the image in $A$ of $\sum_{i=1}^{n-1}\, \pm\qdot
A(in)X_{in} \in T$. 

Note that $R/\ideal{c}\cong \oqmn/\ideal{D_n, X_{1n}}$.  We show that
this is a domain by showing that Jordan's Theorem applies. Note that $c$
is normal (in fact, central) in $R$.  All we need to do is to observe
that $e$ is regular modulo the ideal $dA = Ad$ of the ring $A$. 
However, $A/dA$ is isomorphic to an iterated Ore extension of
${\cal O}_q(M_{n-1})/\ideal{D_{n-1}}$.  
Now ${\cal O}_q(M_{n-1})/\ideal{D_{n-1}}$ is a
domain by \cite[Theorem 2.5]{GL} or by \cite[Example 2]{J}; and so
$A/dA$ is a domain.  Thus, all we have to do is to show that $e\not\in
dA$.

Suppose that $e\in dA$.  Then $e = \widebar{D_{n-1}\alpha}$ for some
$\alpha \in T$. Taking pre-images, we obtain 
\[
\sum_{i=1}^{n-1}\, \pm\qdot A(in)X_{in} = A(nn)\alpha + \beta X_{1n}
\]
for some $\beta \in T$.  
Now, each term on the left hand side is an element of 
${\cal O}_q(M_n)$ of bidegree $(1,\dots,1;1,\dots,1)$ in the 
${\mathbb Z}^n\times{\mathbb Z}^n$ grading of ${\cal O}_q(M_n)$ 
described in [1, 1.5]. 
Hence, we may assume that each term on the right hand side also
has this bidegree.  This implies that $\alpha = \lambda X_{nn}$ for
some scalar $\lambda \in \k$.  Since $X_{nn} \not\in T$, this implies
that $\alpha = 0$.  
In the resulting equation 
\[
\sum_{i=1}^{n-1}\, \pm\qdot A(in)X_{in} = \beta X_{1n}
\]
each of the terms on the left hand side is a preferred product. However,
if we write $\beta = \sum \lambda_i[T_i|T'_i]$ in terms of the preferred
basis then $\beta X_{1n} = \sum \lambda_i[T_i|T'_i]X_{1n}$ is again in
preferred form.  The equation 
\[
\sum_{i=1}^{n-1}\, \pm\qdot A(in)X_{in} = \sum \lambda_i[T_i|T'_i]X_{1n}
\]
then contradicts the independence of the preferred basis. 

Thus $e\not\in dA$ and we have all of the hypotheses of Jordan's Theorem
and conclude that $R/\ideal{c}\cong \oqmn/\ideal{D_n, X_{1n}}$ is a
domain, as required. \QED

It is interesting to note that when $n=2$ the above result fails.  In
this case, the ideal $\langle \det_q , X_{12}\rangle$ is semiprime; in
fact, it is the intersection of two completely prime ideals, each of
which is fixed by the automorphism determined by $X_{1n}$, so
Lemma~\ref{lemme-mo} is applicable.  However, this case has already been
dealt with in \cite{R}. 

\begin{theorem} $\oqmn/\langle \det_q \rangle$ is a maximal order in its
division ring of quotients for each $n\geq 2$. 
\end{theorem}

\proof When $n=2$ this is proved in \cite[Th\'eor\`eme 2.3.11]{R} (and
can be proved directly from Theorem~\ref{main-thm} 
by the reasoning in the previous paragraph).  
An inductive argument similar to that used in Theorem~\ref{main-thm} 
finishes the proof. \QED

\section{Appendix: some useful relations}

In this section, we collect some useful relations. They are
essentially derived from results of \cite{PW}. 

\begin{relation} \label{2X2-q-matrix}  
\rm If $\left( \begin{array}{cc} a&b \\c&d \end{array} \right)$
is a $2 \times 2$ quantum matrix, then the following relation holds:\\
\[ad-q^{2}da=(1-q^{2})(ad-qbc).
\]
\end{relation}

\begin{relations} \label{3X3-q-matrix}  \rm If 
$\left( 
\begin{array}{ccc} 
X_{11}&X_{12}&X_{13} \\
X_{21}&X_{22}&X_{23} \\
X_{31}&X_{32}&X_{33}
\end{array} 
\right)$ 
is a $3 \times 3$ quantum matrix, then the following relations hold.\\

\noindent
1) For $i=2,3$:\\
1.1) $X_{11}[i,1|1,3]
=[i,1|1,3]
X_{11}$, \hfill (see [PW] (4.6.1)) \\
1.2) $X_{11}[i,1|2,3]
=q[i,1|2,3]
X_{11}$, 
\hfill (see [PW] (4.5.1)(1)) \\
1.3) $X_{12}[i,1|1,3]
-q[i,1|1,3]
X_{12}
=(q^{-1}-q)X_{11}[i,1|2,3]
$, \hfill 
(see [PW] (4.5.1)(1)) \\                                
1.4) $X_{12}[i,1|2,3]
=[i,1|2,3]
X_{12}$.  \hfill (see [PW] (4.6.1)) \\

\noindent
2) For $j=1,2$:\\ 
2.1) $X_{33}[3,1|j,3]
=
[3,1|j,3]
X_{33}$, \hfill (see [PW] (4.6.1)) \\
2.2) $X_{33}[2,1|j,3]
=q^{-1}[2,1|j,3]
X_{33}$, 
\hfill (see [PW] (4.5.1)(3)) \\ 
2.3) $X_{23}[3,1|j,3]
-q^{-1}[3,1|j,3]
X_{23}
=(q-q^{-1})X_{33}[2,1|j,3]
$, 
\hfill (see [PW] (4.5.1)(3)) \\
2.4) $X_{23}[2,1|j,3]
=[2,1|j,3]
X_{23}$.  
\hfill (see [PW] (4.6.1)) \\
\end{relations}





\noindent T. H. Lenagan:
School of Mathematics, University of Edinburgh,
\\ James Clerk Maxwell Building, King's Buildings, Mayfield Road,
\\Edinburgh EH9 3JZ, Scotland
\\E-mail: tom@maths.ed.ac.uk
\\
\\
L. Rigal: Universit\'e Jean Monnet
(Saint-\'Etienne), Facult\'e des Sciences et
\\Techniques, D\'e\-par\-te\-ment de Math\'ematiques
\\ 23 rue du Docteur Paul Michelon
\\42023 Saint-\'Etienne C\'edex 2,\\France
\\E-mail: Laurent.Rigal@univ-st-etienne.fr

\end{document}